\documentclass[leqno,12pt, dvipdfmx]{amsart} 

\usepackage{amssymb}
\usepackage{amsmath}
\usepackage{amsthm}
\usepackage[pdftex]{graphicx}
\usepackage[all]{xy}
\usepackage{mathtools}
\usepackage{ascmac}
\usepackage{mathrsfs}
\usepackage{float}
\usepackage{pdfpages}
\usepackage{here}
\usepackage[top=20truemm, bottom=20truemm, left=30truemm, right=30truemm]{geometry}
\theoremstyle{plain} 
\newtheorem{theorem}{\indent\sc Theorem}[section]
\newtheorem{lemma}[theorem]{\indent\sc Lemma}

\newtheorem{proposition}[theorem]{\indent\sc Proposition}
\newtheorem{claim}[theorem]{\indent\sc Claim}

\theoremstyle{definition} 
\newtheorem{definition}[theorem]{\indent\sc Definition}

\newtheorem{observation}[theorem]{\indent\sc Observation}

\newtheorem{question}[theorem]{\indent\sc Question}
\renewcommand{\proofname}{\indent\sc Proof.}

\newcommand{\del}{\partial}
\newcommand{\delbar}{\overline{\partial}}

\newcommand{\U}{\mathrm{U}(1)}
\newcommand{\twopii}{2\pi\sqrt{-1}}
\newcommand{\fractwopii}{\frac{1}{2\pi\sqrt{-1}}}
\newcommand{\Diff}{\mathrm{Diff}^\omega_+(S^1)}
\newcommand{\R}{\mathbb{R}}
\newcommand{\modZ}{\ (\mathrm{mod}\ \mathbb{Z})}
\newcommand{\modtwopiiZ}{\ (\mathrm{mod}\ 2\pi\sqrt{-1}\, \mathbb{Z})}
\newcommand{\one}{\textup{\mbox{1}\hspace{-0.25em}\mbox{l}}}
\newcommand{\fhat}{\widehat{f}}
\newcommand{\psihat}{\widehat{\psi}}

\let\Re\relax
\DeclareMathOperator{\Re}{Re}
\let\Im\relax
\DeclareMathOperator{\Im}{Im}

\begin{document}

\title[Linearization of transiton functions]
{Linearization of transition functions along a certain class of Levi-flat hypersurfaces}

\author[S. Ogawa]{Satoshi Ogawa} 
%
\address{
Department of Mathematics, Graduate School of Science, Osaka Metropolitan University \endgraf
3-3-138, Sugimoto, Sumiyoshi-ku Osaka, 558-8585 \endgraf
Japan
}
\email{sn22894n@st.omu.ac.jp}

\maketitle
\begin{abstract}
We pose a normal form of transition functions along some Levi-flat hypersurfaces obtained by suspension. 
By focusing on methods in circle dynamics and linearization theorems, we give a sufficient condition to obtain a normal form as a geometrical analogue of Arnol'd's linearization theorem. 
\end{abstract}
\section{Introduction}
We study a neighborhood of a {\em Levi-flat hypersurface}. 
Let $X$ be a non-singular complex surface. 
We say a real hypersurface $M$ of $X$ is {\em Levi-flat} if and only if {the Levi-form of $M$ }vanishes identically. 
Especially, $M$ is Levi-flat when $M$ admits a system of local defining functions $\{\rho_j\}_j$ such that each $\rho_j$ is  pluriharmonic  (i.e.  $\del\delbar\rho_j = 0$). 

One of our interests is {the} {\it linearization} of transition functions on a neighborhood of a Levi-flat hypersurface. 
In the study of 1-dimensional complex dynamical systems, it is important to consider whether or not one can find a coordinate by which a function can be regarded as a linear map at a neighborhood of a given fixed point or an invariant curve (linearization, see in \S2). 
In this paper, we will find a normal form of the complex structure of a neighborhood of a Levi-flat hypersurface $M$ by applying a technique for linearization around a circle to the transition functions of a coordinate functions on a neighborhood of $M$. 
In order to apply such a technique for {the} linearization around a circle, we will focus on a certain class of Levi-flat hypersurfaces, which are constructed by {\em suspension construction}. 

Let $Y$ be a non-singular compact complex curve and $\Diff$ be the group of orientation preserving $\mathcal{C}^\omega$-diffeomorphisms of $S^1$, where $S^1$ is the unit circle $\{z \in \mathbb{C} \mid |z| = 1\}$.  
For a given action of the fundamental group $\kappa\colon \pi_1(Y, \ast) \to \Diff$, we consider the quotient space $M$ defined by $Y_{\mathrm{univ}} \times S^1 / \sim$, where $\sim$ is the relation induced from the action $\kappa$  (i.e. $(z,x) \sim (z \cdot \gamma, \kappa(\gamma)(x))$ for $\gamma \in \pi_1(Y, \ast)$) .
Then $M$ is said to be obtained by {\em suspension construction} of $\kappa\colon \pi_1(Y, \ast) \to \Diff$. 

Assume that $M$ is embedded into a non-singular complex surface $X$. 
Let $\mathcal{U} = \{U_j\}$ be a finite covering of $Y$ and $\pi\colon M\to Y$ be the projection. 
For technical reasons, we assume that there exists a holomorphic submersion $P\colon V \to Y$ on a neighborhood $V$ of $M$ in $X$ which satisfies $P|_M = \pi$.

	\begin{definition}\label{GSLF}
We say that $\{(V_j, (z_j, w_j))\}$ is {\em a good system of local functions of width} {$\sigma>0 $  } if and only if it satisfies the following conditions. \\
$(i)$ {For each $U_j$, a local coordinate $z_j$ of $U_j$ and $M_j\coloneqq \pi^{-1}(U_j)$, let $V_j$ be a neighborhood of $M_j$ in $X$ which satisfies $\bigcup_j V_j \subset V$, and $(z_j, w_j)$ be a local coordinate on $V_j$, } {where the coordinate $(z_j, w_j)$ on $V_j$ is given by pullback of the local coordinate of $U_j$ by $P$. }\\
$(ii)$ For each $j$ and $k$, $V_j \cap M = M_j$ holds and {$V_j \cap V_k  = \emptyset$ holds if $M_j \cap M_k = \emptyset $}.  \\
$(iii)$ There exists a positive number $\sigma_j \geq \sigma$ which satisfies the following condition for each $j$: There exists a biholmorphism from $V_j$ to $U_j \times \{e^{-\sigma_j} < |w_j| < e^{\sigma_j}\}$ which makes the following diagram commutative.
\begin{figure}[H]
\centering
\includegraphics[width =6.5cm]{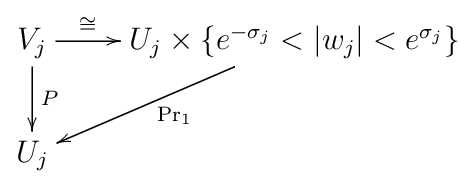}
\end{figure}
$(iv)$ For each $j$, $M_j =\{(z_j, w_j) \in V_j \mid |w_j| = 1\}$ holds. 
	\end{definition}
In what follows, we always assume that a system $\{Y, M, \mathcal{U}, X, \pi, P\}$ has a good system $\{V_j, (z_j, w_j)\}$ of local functions of width $\sigma$. 
Then, the hypersurface $M$ has the local defining function determined by $\log|w_j|$, from which it follows that $M$ is a Levi-flat hypersurface of $X$. 
We recall that $M$ has a structure of $S^1$-bundle over $Y$. 
	We will say that a system $\{Y, \mathcal{U}, \kappa, X, P, \{(V_j (z_j, w_j))\}, \sigma\}$ is  {\em linearizable} if there exists a good system of local function system $\{(V'_j, (z_j, w'_j))\}$ of width $\sigma'$ for $\{Y, \mathcal{U}, \kappa, X, P\}$ such that the transition function is written as $w'_k = t_{jk}w'_j$ on each $V'_{jk} \coloneqq  V'_j \cap V'_k$, where $t_{jk} \in \U$. 
	
Let 
\[
	\begin{cases}
		{z_k = z_k(z_j)}\\
		{w_k = f_{kj}(w_j)}
	\end{cases}
\]
be the transition on $V_{jk} \coloneqq V_j \cap V_k$.
Note that  $z_k$ does not depend on $w_j$ by the condition (iii) in Definition \ref{GSLF}. 
Note also that the transition function $f_{kj}$ does not depend on $z_j$ (see {Lemma \ref{bkj}}). 
We call $f_{jk}$ a {\em transversal transition function} of $V_{jk}$ for a good systems of local functions $\{(V_{j}, (z_j, w_j))\}$. 
The funcition {$f_{kj}|_{S^1}$} is an element of $\Diff$ with a variable $w_j$, where we regard $S^1$ as $\{|w_j| =1\}$. 
Our aim in this paper is to investigate the linearization of transition functions $\{f_{kj}\}_{j, k}$. 

Let $\alpha_{kj}$ and $b_{kj|n}$ be coefficients of the expansion
\[
\log\frac{f_{kj}(w_j)}{w_j} \equiv  \alpha_{kj} + \sum_{n \neq 0}b_{kj|n}w_j^n \modtwopiiZ. 
\]
Note that $\alpha_{kj} \in \sqrt{-1}\mathbb{R}$, since $f_{kj} \in \Diff$ (see the argument in the proof of Lemma \ref{S1}).  
If $\alpha(\{(V_j, (z_j, w_j)\}) \coloneqq \{(U_{jk}, e^{\alpha_{kj}})\} \in \check{C}^1(\mathcal{U}, \U)$ satisfies the 1-cocycle condition, $N = [\{(U_{jk}, e^{\alpha_{kj}})\}] \in \check{H}^1(\mathcal{U}, \U)$ can be regarded as a unitary flat line bundle over $Y$, where $\U = \{t \in \mathbb{C} \mid |t| =1\}$.  
We denote by $b_{kj|n}(\{(V_j, (z_j, w_j))\})$ the non-zero order coefficients $b_{kj|n}$ for a good system of local functions $\{(V_j, (z_j, w_j))\}$.  

The main result is the following. 
	\begin{theorem}\label{main_thm}
Let $Y$ be a compact complex curve, $\mathcal{U} = \{U_j\}$ a finite open covering of $Y$, and $\pi \colon M \to Y$ an $S^1$-bundle over $Y$ constructed by suspension associated to an action $\kappa \colon \pi_1(Y, \ast) \to \Diff$. 
Let $X$ be a complex surface which has $M$ as a Levi-flat hypersurface. 
Assume that there exists a holomorphic submersion $P \colon  V \to Y$ which satisfies $P|_M = \pi$, where $V$ is a neighborhood of $M$ in $X$. 
Assume also that there exists a good system of local functions $\{(V_{j, 0}, (z_j, w_{j, 0}))\}$ of width $\sigma_0$ for $Y, \mathcal{U}, \kappa, X$ and $P$. 
Then, the system $\{Y, \mathcal{U}, \kappa, X, P, \{(V_{j, 0}, (z_j, w_{j, 0}))\}, \sigma_0\}$ is linearizable if the following conditions $(i), (ii), (iii)$ hold. \\
$(i)$ The 1-cochain $\alpha(\{(V_{j, 0}, (z_j, w_{j, 0}))\}) = \{(U_{jk}, e^{\alpha_{kj}})\} \in \check{C}^1(\mathcal{U}, \U)$ satisfies the 1-cocycle condition and $N= [\{(U_{jk}, e^{\alpha_{kj}})\}]$ satisfies $(C_0, \mu, K)$-Diophantine condition, in the sense of Definition \ref{C-dioph} (see \S2), where $C_0 > 0, \mu >1$, and $K$ is constant determined only by $Y$ and $\mathcal{U}$. \\
$(ii)$ For non-zero order coefficients $b_{kj|n, 0} = b_{kj|n}(\{(V_{j, 0}, (z_j, w_{j, 0}))\})$ associated to the transversal transition function $f_{kj, 0}$ of $\{(V_j, (z_j, w_{j, 0}))\}$, there exists a constant $\eta_0 \in (0, \min\{\pi, (1-\mu^{-\frac{1}{\mu+1}})\frac{\sigma_0}{4}\})$ such that 
\[
\max_{j, k} \sup_{e^{-\sigma_0}<|p|<e^{\sigma_0}} \left| \sum_{n\neq0} \, b_{kj|n, 0} \, p^n \right| < \min \left\{ \eta_0, \ \frac{\eta_0^{\mu+1}}{(1+e^{\sigma_0})C_1\mu}\right\}
\] 
holds, where $C_1$ is a constant which depends only on $C_0, \mu$ and $\sigma_0$. \\
$(iii)$ For any good system of local functions $\{(V_j, (z_j, w_j)\}$ which has $[\alpha(\{(V_j, (z_j, w_j))\})] = N$ as a unitary flat line bundle over $Y$, a 1-cochain $\{(U_{jk}, b_{kj|n})\} \in \check{C}^1(\mathcal{U}, \mathcal{O}_Y(N^n))$ satisfies the 1-coboundary condition. \\
	\end{theorem}
Comparing with Arnol'd's linearization theorem (Theorem \ref{Arnold_linearization}) and Ueda's linearization theorem (Theorem \ref{linearization_2}), we will explain the conditions $(i)$, $(ii)$, and $(iii)$. 
The condition $(i)$ is the more detailed version of {the} Diophantine condition. The condition $(ii)$ corresponds to the assumption for the estimate of a perturbation in Arnol'd's linearization theorem. The condition $(iii)$ corresponds to vanishing of obstruction classes in Ueda's proof. 

Main result can be applied to finding a criterion of simultaneous linearization of circle diffeomorphisms (see \S\ref{example}). In this sense, Theorem \ref{main_thm} can be regarded as a generalization of Arnol'd's linearization theorem.  

{
In \cite{KU}, Koike and Uehara constructed Levi-flat in K3 surfaces and showed there exist a foliated small neighborhood of it. 
From the view point of Koike and Uehara's result, this result also can be regarded as one of a result for a foliated neighborhood along Levi-flat hypersurfaces. 
}

Our idea and main result can be explained as a geometrical analogue of Arnold's linearization theorem. 
For proving Theorem \ref{main_thm}, we use Kolmogorov-Arnol'd-Moser (KAM) theory, which is used in the proof of Arnold's linearization theorem (Theorem \ref{Arnold_linearization}) in \cite{CG} \cite{SM}.  
In \cite{U}, Ueda investigated linearization on a neighborhood of a compact complex curve embedded holomorphically in a complex surface as a geometrical analogue of Siegel's linearization theorem (Theorem \ref{Siegel_linearization}). 

In \S\ref{preliminaries}, we introduce preliminaries about linearization theorems. 
In \S\ref{circle_dynamics} we will explain the expansion of transition functions. 
In \S\ref{linearization_1} and \S\ref{linearization_2}, we will see two linearization theorems in one-dimensional dynamics and Ueda's linearization theorem. 
In \S\ref{proof}, we will apply a method in \S\ref{linearization_1} to a Levi-flat hypersurface constructed by suspension and show Theorem \ref{main_thm}. 	
In \S\ref{example}, I will introduce a simple example on the main result and obtain a sufficient condition for simultaneous linearization of circle diffeomorphisms. 
In \S\ref{discussion}, I will discuss the relation between the expansion of transition functions and a rotation number. 	
\section{Preliminaries}\label{preliminaries}

\subsection{Dynamical system of a circle diffeomorphism}\label{circle_dynamics}
In this section, we will review some fundamental facts on 1-dimensional dynamical systems.   

For $f \in \Diff$, we say that a homeomorphism $F \colon \R \to \R$ is a {\em lift} of $f$ if $F$ satisfies a relation $e^{\twopii F(x)} = f(e^{\twopii x})$. 
The following limit is called the {\em rotation number} of $f$:
\[
\rho(f) \coloneqq \lim_{m\to\infty} \frac{F^{\circ m}(x) -x}{m} \modZ.
\]
It is known that $\rho(f)$ exists and is independent of a choice of $F$ and a point $x \in \R$. 
We have the Fourier expansion of the lift of $f \in \Diff$ as below. 
\[
F(x) = x +F_0 + \sum_{n \neq 0} F_n e^{\twopii nx}. 
\]
By letting $w = e^{\twopii x}$ and taking the logarithm, we obtain 
\[
\log \frac{f(w)}{w} = \twopii F_0 + \sum_{n \neq 0}\twopii F_n w^n. 
\]
Note that there exists a branch of $\displaystyle{\log\frac{f(w)}{w}}$ globally on the annulus $A \coloneqq \{e^{-\sigma} < |w| < e^{\sigma}\}$ for $\sigma > 0$. 
	\begin{lemma}\label{expansion}
There exists $h \in H^0(A, \mathcal{O}_A)$ such that $e^h = g$ holds, where $\displaystyle{g(w) = \frac{f(w)}{w}}$. 
	\end{lemma}
	\begin{proof}
By considering the exponential sheaf exact sequence
\begin{figure}[H]
\centering
\includegraphics[width = 8.0cm]{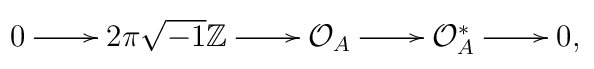}
\end{figure}
one obtains the exact sequence of cohomology groups
\begin{figure}[H]
\centering
\includegraphics[width = 9.75cm]{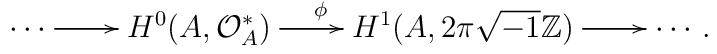}
\end{figure}
It is sufficient to check $\phi(g) = 0$. 
We can calculate $\phi(g)$ as
\[
\phi(g) = \int_\gamma \frac{g'(w)}{{g(w)}} dw, 
\]
where $\gamma$  is a loop in the annulus $A$ which generates the fundamental group $\pi_1(A, \ast)$.   
We can easily check 
\[
\frac{g'(w)}{g(w)} = \frac{f'(w)}{f(w)} - \frac{1}{w}.
\]
Since $f \in \Diff$,  it is shown that $\phi(g) = 0$. 
	\end{proof}
\subsection{Linearization theorems in 1-dimensional dynamics}\label{linearization_1}
Here we survey some studies of the local behavior of a holomorphic function on a neighborhood of a fixed point. 
Let $f$ be a holomorphic function which admits the origin as a fixed point. 
Suppose that $f$ has the expansion $f(z) =\Lambda z + b_2 z^2 + b_3 z^3 +\cdots$ on a neighborhood of the origin. 
A number $\Lambda$ is called the {\em multiplier} of $f$ at the fixed point.
It is known that a classification of a fixed point is given according to the multiplier $\Lambda$.  
	\begin{definition}
An irrational number $\theta$ is said to be {\em Diophantine} if and only if there exist $c >0$ and $B >1$ so that
\[
\left| \theta - \frac{p}{q} \right| \geq \frac{c}{q^B}
\] 
for any rational number $p/q \ (p, q \in \mathbb{Z}, q > 0)$. 
	\end{definition}
We assume that $\Lambda = e^{\twopii \theta}$ holds for an irrational number $\theta$. 
In this case, the fixed point 0 is called an {\em {irrationally neutral} fixed point}. The following theorem is known as an important linearization  theorem at an {irrationally} neutral fixed point.  
	\begin{theorem}[Siegel's linearization theorem \cite{S}]\label{Siegel_linearization}
Let $f$ be a holomorphic function which has the origin as an irrational fixed point with multiplier $e^{\twopii \theta}$. 
	{If} $\theta$ satisfies {the} Diophantine condition, then there exists a holomorphic map $\psi$ on a neighborhood of the origin such that $\psi$ satisfies $\psi(0) = 0$, $\psi'(0) = 1$, and $(\psi^{-1}\circ  f \circ \psi)(z) = e^{\twopii \theta}z$. 
	\end{theorem}
When $\psi$ as in Theorem \ref{Siegel_linearization} exists, we say that $f$ is {\em linearizable} on a neighborhood of 0. To show the linearizablity of $f$ is equivalent to solve the following equation called {\em Schr\"{o}der's equation}: 
\[
\psi( e^{\twopii \theta} z  ) = f ( \psi(z) ).
\]
We define $\fhat$ and $\psihat$ by $f(z) = e^{\twopii \theta}z + \fhat(z)$ and $\psi(z) = z + \psihat(z)$. 
Furthemore suppose that the function $\psihat$ can be written as $\psihat(z) = a_2z^2 + \cdots$. Then, Schr\"{o}der's equation can be rewritten as {$\psihat( e^{\twopii \theta}z) -  e^{\twopii \theta}\psihat(z) = \fhat(\psi(z))$}, which leads that $a_n$ can be determined by $a_2, \cdots a_{n-1}, b_2, \cdots, b_n$ inductively. 
In Siegel's original method, he estimated $a_n$ and proved the convergence of $\psihat$. 

The following Theorem \ref{Arnold_linearization} is a counterpart of {Theorem} \ref{Siegel_linearization} in circle dynamics. 
	\begin{theorem}[Arnold's linearization theorem {\cite{A}},  Theorem 7.2 of \S2.7 in {\cite{CG}},  cf. Theorem 12.3.1 {\cite{KH}}]\label{Arnold_linearization}
Let $\alpha \in [0,1)$ be a number which satisfies {the} Diophantine condition and $\sigma$ be a positive constant.  Then there exists a positive constant $\delta$ such that, if $f$ is any element of $\Diff$ with $\rho(f) = \alpha$ which extends to be analytic and univalent on the annulus $\{e^{-\sigma} < |z| < e^\sigma\}$ and satisfies $|f(z) - e^{\twopii\alpha}z| < \delta$ on $\{e^{-\sigma} < |z| < e^\sigma\}$, then $f$ is linearizable on the annulus $\{e^{-\sigma'} < |z| < e^{\sigma'}\}$, where $0< \sigma' < \sigma$. 
	\end{theorem}
Arnol'd's theorem can be regarded as the linearizaition along the unit circle.
This theorem can be proven by a different strategy from that of Siegel for Theorem \ref{Siegel_linearization}. 
The proof in \cite{CG} is based on a simple case in KAM theory. 
In Arnol'd's proof, we inductively change the coordinates along the unit circle and estimate the non-linear part of $f$. 
Details of this technique will be explained in \S3 ({\em{cf}}: \cite{H}). 

\subsection{Ueda's linearization theorem}\label{linearization_2}
Siegel's linearization theorem can be generalized in a geometric sense, which is known as Ueda's linearization theorem (\cite{U}). 
Let $C$ be a compact complex curve which is holomorphically embedded in a complex surface $S$ with the topologically trivial normal bundle $N_{C/S}$. 
Note that $\{S_j\}$ is a finite open covering of a neighborhood of $C$ in $S$. 
For $S_j$, let $s_j$ be a defining function of $S_j \cap C $ in $S_j$. 
We suppose that, for any $j$ and $k$,  there exists $t_{kj} \in \U$ such that $s_k = t_{kj} s_j + O(s_j^2)$. 
Ueda gave a sufficient condition for the existence of an open covering $\{S'_j\}$ and a system of defining functions  $\{s'_j\}$ such that $s'_k = t_{kj}s'_j$ holds by using {the} Diophantine condition of the normal bundle $N_{C/S}$. 
{The} Diophantine condition {in the sence of Ueda} is defined by focusing on an invariant distance $d$ of $\mathrm{Pic}^0(C)$, where the invariant property of the distance $d$ means that the following holds for any  $E_1, E_2, G \in \mathrm{Pic}^0(C)$:
\[
d(E_1, E_2) = d(E_1^{-1}, E_2^{-1}) = d(E_1 \otimes G, E_2 \otimes G).
\]

	\begin{definition}
For $E \in \mathrm{Pic}^0(C)$ which satisfies $E^{\otimes n} \neq \one$ for any $n$, $E$ is said to be {\em Diophantine} if and only if the following holds:
\[
-\log d(\one, E^{\otimes n} ) = O(\log n) \ (n \to \infty).
\]
	\end{definition}
By using {the} Diophantine condition of a flat line bundle and the following theorem, he found some estimates of coefficients of transition functions at a neighborhood on $C$ and showed the linearization theorem (\cite{U}).

	\begin{theorem}[Lemma 4 in \cite{U}]\label{U}
	{Let $\mathcal{W} = \{W_j\}$ be a finite open covering of $C$. }
There exists a positive constant $K = K(C, \mathcal{W})$ such that the following holds for any flat line bundle $E$ and any $\mathcal{G} \in C^0(\mathcal{W}, \mathcal{O}(E))$: 
\[
d(\one, E)||\mathcal{G}|| \leq K ||\delta\mathcal{G}||, 
\] 
where $\delta$ is the coboundary map from $\check{C}^0(\mathcal{W}, \mathcal{O}(E))$ to $\check{C}^1(\mathcal{W}, \mathcal{O}(E))$ and the norms is defined by  
\[
||\mathcal{G}^0|| \coloneqq \max_j \sup_{p\in W_j} |g_j(p)| 
\]
and
\[
 ||\mathcal{G}^1|| \coloneqq \max_{j, k}\sup_{p \in W_j \cap W_k}|g_{jk}(p)|
\] 
for a 0-cochain $\mathcal{G}^0 = \{(W_j, g_j)\} \in C^0(\mathcal{W}, \mathcal{O}(E))$ and a 1-cochain $\mathcal{G}^1 = \{(W_j \cap W_k,  g_{jk})\} \in C^1(\mathcal{W}, \mathcal{O}(E))$. 
	\end{theorem}
In this paper, by using $K$ in Theorem \ref{U}, we will classify the Diophantine condition in more detail as follows. 
	\begin{definition}\label{C-dioph}
Let $C$ be a compact complex curve, $\mathcal{W}$ be a finite open covering of $C$, and $K$ be the constant as in Theorem \ref{U}. 
The unitary flat bundle $E$ on $C$ is said to satisfy $(C_0, \mu, K)$-Diophantine condition for $C_0 > 0$ and $\mu > 1 $ if and only if
\[
C_0 n^{\mu-1}d(\one, E^{\otimes n}) \geq K
\] 
holds for any $n = 1, 2, \dots$. 
	\end{definition}
%
%
\section{Proof of main theorem}\label{proof}
\subsection{Outline of proof}
In this section, for $Y, \mathcal{U}, \kappa$, $X$, and $P$ given in \S1, we will prove Theorem \ref{main_thm}. 
From Theorem \ref{U}, we obtain a constant $K = K(Y, \mathcal{U})$. 

Let $\{(V_{j, 0}, (z_j, w_{j, 0}))\}$ be an initial good system of local functions of width $\sigma_0 > 0$ over $Y, \mathcal{U}, \kappa$, $X$, and $P$. 
Assume that  $\{Y, \mathcal{U}, \kappa, X, P, \{(V_{j, 0}, (z_j, w_{j, 0}))\}, \sigma_0\}$ satisfies  the  assumption $(i), (ii), (iii)$ in Theorem \ref{main_thm}. 
In \S\ref{psi}, we will explain how to retake of coordinates so that the $\mathcal{L}^\infty$-norm of the non-linear part of the transversal transition function becomes smaller. 
By making width $\sigma_0$ slightly smaller, together with the assumption $(iii)$ in Theorem \ref{main_thm}, one can obtain a function which renew a coordinate $w_{j, 0} \mapsto w_{j, 1}$ on a neighborhood of $M_j$ not by changing a unitary flat line bundle $N = [\{(U_{jk}, {e^{\alpha_{kj}}})\}]$ over $Y$. 
By an inductive procedure, we obtain $\{(V_{j, m}, (z_j, w_{j, m}))\}$ as a good system of local functions of width $\sigma_m$ retaken $m$-times from the initial system $\{(V_{j, 0}, (z_j, w_{j,0}))\}$. 
Since $N$ is not changed by the procedure, a transversal transition $f_{kj, m} \in \Diff$ on $V_{kj, m} =V_{j, m} \cap V_{k, m}$ has the same non-zero order part of $f_{kj, 0} $as $\alpha_{kj}$. 
Therefore, by Lemma \ref{expansion}, it is checked that $f_{kj, m}$ has the Laurent expansion
\[
\log \frac{f_{kj ,m}(w_{j, m})}{w_{j, m}} = \alpha_{kj} + \sum_{n\neq 0}b_{kj|n, m}w_{j, m}^n \modtwopiiZ.
\]
We will denote the sum $\displaystyle{\sum_{n\neq 0}b_{kj|n, m}w_{j, m}^n}$ by $\fhat_{kj, m}$. 
We define the norm by
\[
||\fhat_{kj, m}||_{\sigma_m} \coloneqq \sup_{e^{-\sigma_m} < |w_{j, m}| < e^{\sigma_m}} \left|\fhat_{kj, m}(w_{j, m}) \right|. 
\]
Our goal is to prove that $||\fhat_{kj, m}||_{\sigma_m}$ converges to zero. 
For proving the main theorem, it is sufficient to show the following statement. 
	\begin{theorem}\label{main_theorem2}
Assume that $\{Y, \mathcal{U}, \kappa, X, P, \{(V_{j, 0}, (z_j, w_{j, 0}))\}, \sigma_0\}$ satisfies the assumption in Theorem  \ref{main_thm}. 
Define $\delta_m, \sigma_m$ and $\eta_m$ inductively by
		\begin{gather*}
\delta_0 = \min \left\{ \eta_0, \ \frac{\eta_0^{\mu + 1}}{(1 + e^{\sigma_0})C_1 \mu}\right\}, \\
\eta_{m+1} = \mu^{-\frac{1}{\mu+1}} \eta_m, \\
\delta_{m+1} = (1 + e^{\sigma_0})C_1 \frac{\delta_m^2}{\eta_m^{\mu+1}}, 
		\end{gather*}
and
\[
\sigma_{m+1} = \sigma_m - 4\eta_m,
\]
where $C_0 >0$ and $\mu > 1$ are constants such that $N$ is $(C_0, \mu, K)$-Diophantine, ${\sigma_0 } > 0$ is the initial width, $\eta_0$ satisfies $0 < \eta_0 < \min \{\pi, (1-\mu^{-\frac{1}{\mu+1}})\frac{\sigma_0}{4}\}$ and $C_1$ is a {constant} which depends only on $C_0, \mu$ and $\sigma_0$ (see Lemma \ref{estimate_psi}) . 

Then, the following holds for any $m$:
		\begin{align*} 
\max_{j, k} ||\fhat_{kj, m}||_{\sigma_m} < \delta_m \leq \min \left\{ \eta_m, \ \frac{\eta_m^{\mu+1}}{(1+e^{\sigma_0})C_1 \mu } \right\}. 
		\end{align*}
	\end{theorem}
Recall that $K = K(Y, \mathcal{U})$ is invariant under the inductive procedure. 
Note that $\delta_m \to 0$, since $\eta_m \to 0$ and $\delta_m < \eta_m$. 
Vanishing of the limit of $\delta_m$ allows to deduce that $||\fhat_{kj, m}||_{\sigma_m} \to 0$ as $m \to \infty$. 
From the definition of $\{\sigma_m\}$, one can check directly that the limit of $\sigma_m$ is a positive constant. 

In what follows, $\{(V_{j, m}, (z_j, w_{j, m}))\}$ is a good system of local functions retaken $m$-times and satisfies the inductive assertion as above. 
In \S\ref{f}, we will give some estimates of transition functions $f_{kj, m}$. 
In \S\ref{psi}, we will explain a function of retaking coordinates $w_{j, m} \to w_{j, m+1}$. 
In \S\ref{new}, we will define renewed transition functions $f_{kj, m+1}$ from $f_{kj, m}$ and give an estimate of $f_{kj, m+1}$. 
%
%
\subsection{Review of transition function}\label{f}
The transition on $V_{jk, m} \coloneqq V_{j, m} \cap V_{k, m}$ is given by
\[
	\begin{cases}
		{z_k = z_k(z_j)}\\
		{w_{k, m} = f_{kj, m}(w_{j, m})}
	\end{cases}. 
\]
We denote an expansion of $f_{kj, m} \in\Diff$ by
\[
\log\frac{w_{k, m}}{w_{j, m}} = \log\frac{f_{kj, m}(w_{j, m})}{w_{j, m}} \equiv \alpha_{kj} + \sum_{n\neq 0} b_{kj|n, m}w_{j, m}^n \modtwopiiZ. 
\]

{
Firstly, we check that the transition function $f_{kj, m}$ does not depend on $z_j$. 
In the following lemma, we denote an expansion of $f_{kj, m} \in \Diff$ by
\[
\log\frac{f_{kj, m}(w_{j, m}, z_j)}{w_{j, m}} \equiv \alpha_{kj} + \sum_{n\neq 0} b_{kj|n, m}(z_j)w_{j, m}^n \modtwopiiZ. 
\]
\begin{lemma}\label{bkj}
For any $n \neq 0$, $\overline{b_{kj|n, m}}(z_j) = -b_{kj|-n, m}(z_j)$ holds. 
Especially, the transition function $f_{kj, m}$ does not depend on $z_j$. 
\end{lemma}
\begin{proof}
By considering a contour integral over the unit circle, one has
\[
b_{kj|n, m} (z_j)= \int_0^{1}  \log \left( \frac{f_{kj, m}(e^{\twopii \theta}, z_j)}{e^{\twopii\theta}} \right) e^{-2\pi n\sqrt{-1}\theta}d\theta. 
\]
Taking conjugation of $b_{kj|n, m}$, 
\[
\overline{b_{kj|n, m}}(z_j) = \int_0^{1}  \overline{ \log \left( \frac{f_{kj, m}(e^{\twopii \theta}, z_j)}{e^{\twopii\theta}} \right) } e^{2\pi n\sqrt{-1}\theta} d\theta. 
\]
Since $f_{kj, m}\in \Diff$, $\overline{b_{kj|n, m}(z_j)} = -b_{kj|-n, m}(z_j)$ holds. 
Therefore, $b_{kj|n, m}$ is holomorphic and antiholomorphic function of $z_j$. 
This implies that $b_{kj,m}$ is a constant with respect to $z_j$. 
\end{proof}
}

	\begin{lemma}\label{estimate_b} For any $n \neq 0$, the following holds: 
\[
|b_{kj|n, m}| \leq ||\fhat_{kj}||_{\sigma_m} e^{-|n|\sigma_m}.
\]
	\end{lemma}
{\proofname} From the definition,
\[
{
|b_{kj|n, m}| \leq \frac{1}{2\pi} \int_\gamma \left|\log\frac{f_{kj, m}(w_{j, m})}{w_{j, m}} - \alpha_{kj} \right| \frac{d|w_{j, m}|}{|w_{j, m}|^{n+1}},
}
\]
where $\gamma$ is a generating loop. 
Considering the loop $\gamma = \{|w_{j, m}| = e^{\pm\sigma_m}\}$, we obtain the inequality. \qed \\
\subsection{A function of retaking coordinates}\label{psi}
In this section, we will define the renewed coordinate $\{ (z_j, w_{j, m+1})\}$ from the coordinate $\{(z_j, w_{j, m})\}$ by suitably constructing the function $\psi_{j, m}$: 
\[
w_{j, m} = \psi_{j, m} (w_{j, m+1}). 
\]
The function $\psi_{j, m}$ will be constructed by using the function
\[
\psihat_{j, m}(w' ) = \sum_{n\neq0} a_{j|n, m} w'^n, 
\]
where $a_{j|n,m}$ are suitably chosen constants {to satisfy
\[
\psihat_{j, m} (w') = \log \frac{\psi_{j, m}(w')}{w'}. 
\]
}
Let us explain how to construct of $a_{j|n, m}$.
We obtain $\psi_{j, m}$ from the {\em simplified {Schr\"{o}der's equation}}
	\begin{equation}\label{Schroder_eq}
\fhat_{kj, m}(w_{j, m+1}) + \psihat_{j, m}(w_{j, m+1}) - \psihat_{k, m}(e^{\alpha_{kj}}w_{j, m+1}) = 0
	\end{equation}
\begin{observation}\label{obs}
Let us explain the simplified Schr\"{o}der's equation (\ref{Schroder_eq}).
For simplicity assume that transitions on $V_{jk, m+1}$ is linear: {$w_{k, m+1} = e^{\alpha_{kj}} w_{j, m+1}$}. 
Then, functions of retaking coordinates satisfy 
	\begin{align*}
\alpha_{kj} = \log\frac{w_{k, m+1}}{w_{j, m+1}} &= \log \left( \frac{w_{k, m+1}}{w_{k, m}}\cdot\frac{w_{k, m}}{w_{j, m}}\cdot\frac{w_{j, m}}{w_{j, m+1}} \right) \\
&=  -\psihat_{k, m}(w_{k, m+1}) + \log\frac{w_{k, m}}{w_{j, m}} + \psihat_{j, m}(w_{j, m+1}).
	\end{align*}
Thus, we obtain $\fhat_{kj, m}(w_{j, m}) + \psihat_{j, m}(w_{j, m+1}) - \psihat_{k, m}(w_{k, m+1}) = 0$ as {\em Schr\"{o}der's equation}.
All we have to do is to find the solution of this, which is not easy. 
Instead of solving Schr\"{o}der's equation, we consider (\ref{Schroder_eq}) as a {\em simplified Schr\"{o}der's equation} (replace $w_{j,m}$ with $w_{j, m+1}$ and $w_{k, m+1}$ with $e^{\alpha_{kj}}w_{j, m+1}$). This idea came from the KAM theoretical proof of Theorem \ref{Arnold_linearization}. \qed \\
\end{observation}

By using power series, {the simplified Schr\"{o}der's equation (\ref{Schroder_eq})} turns out that $a_{j|n, m}$ should satisfy
\[
b_{kj|n, m}+a_{j|n, m} - (e^{\alpha_{kj}})^n a_{k|n, m} = 0. 
\]
It is easily checked that this condition is equivalent to the existence of $\{a_{j|n, m}\}$ which satisfies $\delta_{n, m}^0( \{(U_j, a_{j|n, m})\} ) =  \{(U_{jk}, b_{kj|n, m})\}$, where
\[
\delta_{n, m}^0: \check{C}^0 (\mathcal{U}, \mathcal{O}_Y(N^n)) \to \check{C}^1 (\mathcal{U}, \mathcal{O}_Y(N^n))\\
\]
is the coboundary map.
From the assumption $(iii)$ in Theorem\ref{main_thm}, one can find $\{a_{j|n, m}\}$ which satisfies the condition above. 
In this situation, since the compactness of $Y$ and the unitary-flatness of $N$, we can apply Theorem \ref{U} to $\{a_{j|n, m}\}$ and $\{b_{kj|n, m}\}$ to conclude that
\[
d(\one, N^{\otimes n})\max_{j}|a_{j|n, m}| \leq K \max_{j, k}|b_{kj|n, m}|. 
\]
Recall that $K$ does not depend on $n$ and $m$. 
From an invariant property $d(\one , N^{\otimes n}) = d(\one, N^{\otimes (-n)})$, 
\[
C_0 |n|^{\mu-1} d(\one, N^{\otimes n}) \geq K 
\] 
holds for any $n \neq 0$ if $N$ satisfies $(C_0, \mu, K)$-Diophantine condition. 
By combining these estimates, we obtain the following: 
	\begin{lemma}\label{Dioph}
If $N$ satisfies $(C_0, \mu, K)$-Diophantine condition, 
\[
\max_{j}|a_{j|n, m}| \leq C_0|n|^{\mu-1} \max_{j, k} |b_{kj|n, m}|
\]
holds for any $n \neq 0$. 
	\end{lemma}
In this manner, we obtain a function of retaking coordinates 
\[
\psihat_{j, m} (w_{j, m+1}) = \sum_{n\neq0} a_{j|n, m} w_{j, m+1}^n
\]
 (the convergence of $\psihat_{j, m}$ will be proven later in this section). 
For $\sigma' > 0$, we define the norm $|| \cdot ||_{\sigma'}$ by 
\[
|| \psihat_{j, m} ||_{\sigma'} \coloneqq \sup_{e^{-\sigma'} < |w_{j, m+1}| < e^{\sigma'}} |\psihat_{j, m}(w_{j, m+1})|
\]
Next, we estimate the norm of the function of retaking coordinates. 
	\begin{lemma}\label{estimate_psi}
There exists a constant $C_1$ which depends only on $C_0, \mu$, and $\sigma_0$ such that 
\[
|| \psihat_{j, m} ||_{\sigma_m - \lambda} \leq C_1 \cdot  \max_{j, k} ||\fhat_{kj, m}||_{\sigma_m}  \cdot \lambda^{-\mu}.
\] 
holds for any $\lambda \in (0, \sigma_m)$. 
	\end{lemma}
{\proofname} From {Lemma \ref{estimate_b}} and {Lemma \ref{Dioph}, we have
\begin{align*}
|| \psihat_{j, m}||_{\sigma_m - \lambda} 
&= \sup_{e^{-(\sigma_m - \lambda)} < |w_{j, m+1}| < e^{\sigma_m - \lambda}}\left| \sum_{n \neq 0}a_{j|n, m}\,w_{j, m+1}^n \right| \\
&\leq \sum_{n \neq 0} |a_{j|n, m}|\, e^{(\sigma_m - \lambda)|n|} \\
&\leq \sum_{n \neq 0} C_0 |n|^{\mu-1} \cdot \max_{j, k}|b_{kj|n, m}| \cdot e^{(\sigma_m - \lambda)|n|} \\
&\leq 2C_0 \cdot \max_{j, k}||\fhat_{kj, m}||_{\sigma_m} \cdot \sum_{n\geq 1} n^{\mu-1}e^{-n\lambda}. 
\end{align*}
The sum in the right hand side can be calculated as follows:  
\[
\sum_{n\geq1} n^{\mu -1} e^{-n\lambda} \leq e^{-\lambda}\cdot\sum_{n\geq1}n(n+1)\cdots(n+\mu-2) (e^{-\lambda})^{n-1} = \frac{e^{-\lambda}\cdot (\mu-1)!}{(1-e^{-\lambda})^\mu} \leq \frac{(\mu - 1)!}{(1-e^{-\lambda})^\mu}. 
\]
From the convexity of the function $x \mapsto 1-e^{-x}$, for any $\lambda$ which satisfies $0 < \lambda < \sigma_m <\sigma_0$, the following holds: 
\[
\frac{1-e^{-\lambda}}{\lambda} > \frac{1-e^{-\sigma_0}}{\sigma_0}.
\]
Hence, we obtain
\[
||\psihat_{j, m}||_{\sigma_m - \lambda} \leq 2 C_0 \cdot \max_{j, k}||\fhat_{kj, m}||_{\sigma_m} \cdot \frac{{\sigma_0}^\mu(\mu-1)!}{(1-e^{-\sigma_0})^\mu} \cdot \lambda^{-\mu}.
\]
Letting $\displaystyle{C_1 = \frac{2C_0\,\sigma_0^\mu (\mu-1)!}{(1-e^{-\sigma_0})^\mu}}$, the statement is proven. \qed \\

We shall check the well-definedness of retaking coordinates. 
Let $V_{j, m}(\sigma') \coloneqq \{ e^{-\sigma'} < |w_{j, m}| < e^{\sigma'} \}$ and  $V_{j, m+1}(\sigma') \coloneqq \{ e^{-\sigma'} < |w_{j, m+1}| < e^{\sigma'} \}$ {for a positive real number $\sigma'$}. 
For coordinates $w_{j, m}$ and $w_{j, m+1}$, note that the renewed transversal transition function $f_{kj, m+1}$ is defined by
\[
w_{j, m+1} \mapsto f_{kj, m+1}(w_{j, m+1}) = (\psi_{k, m}^{-1}\circ f_{kj, m} \circ \psi_{j, m})(w_{j, m+1}).
\]
Recall that $\sigma_m - \nu \eta_m$ is positive for any  $\nu\in \{1, 2, 3, 4\}$ by definitions of $\{\sigma_m\}$ and $\{\eta_m\}$ {in Theorem \ref{main_theorem2} }.
	\begin{proposition}\label{well-definedness}
The function $f_{kj, m+1}$ is well-defined as a map from $V_{j, m+1}(\sigma_{m} - 4\eta_{m})$ to $V_{k, m+1}(\sigma_{m} - \eta_{m})$. 
	\end{proposition}
{\proofname} We prove this theorem by checking the following properties.
	\begin{enumerate}
		\item $\psi_{j, m}(V_{j, m+1}(\sigma_m-4\eta_m)) \subset V_{j, m}(\sigma_m-3\eta_m)$. 
		\item $f_{kj, m}(V_{j, m}(\sigma_m-3\eta_m)) \subset V_{k, m}(\sigma_m-2\eta_m)$. 
		{
		\item $\psi_{k, m}^{-1}$ is well-defined on $V_{k, m}(\sigma_m - 2\eta_m)$. 
		\item $ \psi_{k, m}^{-1}( V_{k, m}(\sigma_m -2\eta_m) )\subset V_{k, m+1}(\sigma_m-\eta_m)$. 
}
	\end{enumerate}

{
{\sc Proof of (1) and (2). }}
Note that $\psi_{j,m}$ is well-defined on $V_{j, m+1}(\sigma_m - \lambda)$ for $0<\lambda<\sigma_m$ by Lemma \ref{estimate_psi}. 
From Lemma \ref{estimate_psi}, for $\nu = 1, 2, 3, 4$, we obtain
	\begin{align*}
||\widehat{\psi}_{j, m}||_{\sigma_m - \nu\eta_m} &\leq C_1\cdot(\nu\eta_m)^{-\mu}\cdot\left( \max_{j, k}||\widehat{f}_{kj}||_{\sigma_m} \right) \\
&< C_1 \cdot (\nu\eta_m)^{-\mu} \cdot \frac{\eta_m^{\mu+1}}{(1+e^{\sigma_0})C_1\mu} < \eta_m
	\end{align*}
by the inductive assumption $\displaystyle{\max_{j, k}||\widehat{f}_{kj}||_{\sigma_m} < \frac{\eta_m^{\mu+1}}{(1+e^{\sigma_0})C_1\mu}}$.
Therefore one has
\begin{align*}
|\log|\psi_{j, m}(w_{j, m+1})| - \log|w_{j, m+1}|| &= |\Re\, (\psihat_{j, m} (w_{j, m+1}))| \\
&<|{\psihat_{j, m}(w_{j, m+1})}| \\
&< \eta_m
\end{align*}
on $V_{j, m+1}(\sigma_m - 4\eta_m)$, which proves the assertion (1). \\

The assertion (2) is proven from the following inequality for $w_{j, m} \in V_{j, m}(\sigma_m - 3\eta_m)$ and the inductive assumption: 
\begin{align*}
|\log|f_{kj, m}(w_j)| - \log|w_{j, m}| | &= |\Re\, (\fhat_{kj, m}(w_{j, m}))| \\
&<|\fhat_{kj, m}(w_{j, m})| \\
&<\max_{j, k}||\fhat_{kj, m}||_{\sigma_m} \\ 
&< \delta_m < \eta_m. \qed
\end{align*}

\noindent 
To prove (3) and (4), we shall use the following lemma. 
	\begin{lemma}\label{1+e}
\[
\sup_{|\Re\, \zeta| < \sigma_m - \eta_m} \left| \frac{d}{d\zeta}\psihat_{k, m}(e^\zeta) \right|\leq \frac{1}{1+e^{\sigma_0}}. 
\]
	\end{lemma}
{\proofname} By using power series expression of $\psihat_{k, m}$ and the same argument as in the proof of Lemma $\ref{estimate_psi}$, one has
	\begin{align*}
\sup_{|\Re\, \zeta| < \sigma_m - \eta_m} \left| \frac{d}{d\zeta}\psihat_{k, m}(e^\zeta) \right|&\leq 2C_0 \cdot\max_{j, k}||\fhat_{kj, m}||_{\sigma_m} \cdot \sum_{n \geq 1} n^{\mu}e^{-n\eta_m} \\
&\leq 2C_0 \cdot \max_{j, k}||\fhat_{kj, m}||_{\sigma_m} \cdot \frac{\mu!\, e^{-\eta_m}}{(1-e^{-\eta_m})^{\mu+1}} \\
&= C_1 \mu\cdot \max_{j, k}||\fhat_{kj, m}||_{\sigma_m} \cdot \left( \frac{1-e^{-\sigma_0}}{\sigma_0 (1-e^{-\eta_m})}\right)^\mu \cdot \frac{1}{e^{\eta_m}-1}. 
	\end{align*}
Thus we have
	\begin{align*}
\sup_{|\Re\, \zeta| < \sigma_m - \eta_m} \left| \frac{d}{d\zeta}\psihat_{k, m}(e^\zeta) \right| 
&\leq C_1 \mu \cdot \max_{j, k}||\fhat_{kj, m}||_{\sigma_m} \cdot \frac{1}{\eta_m^\mu(e^{\eta_m}-1)}.
	\end{align*}
Therefore the assertion follows from $\displaystyle{\max_{j, k}||\fhat_{kj, m}||_{\sigma_m} < \frac{\eta_m^{\mu+1}}{(1+e^{\sigma_0})C_1\mu}}$ and the inequality $\displaystyle{\frac{\eta_m}{e^{\eta_m} -1} < 1} $. \qed \\
Furthermore we shall prove the following.
	\begin{lemma}\label{S1}
It follows that $\psi_{k, m} (S^1) = S^1$, where $S^1$ is identified with $\{|w_{k, m+1}| =1 \}$.
	\end{lemma}
{\proofname} 
{From lemma \ref{bkj}, the relation $\overline{b_{kj|n, m}} = -b_{kj|-n, m}$ holds. 
This relation leads
\[
\delta^0_{n, m}(\{(U_j, a_{j|n,m})\}) = \delta^0_{n, m}(\{(U_j, - \overline{a_{j|-n,m}})\}) = \{(U_{jk}, b_{kj|n, m})\}.
\] 
}
One also has $H^0(\mathcal{U}, \mathcal{O}_Y(N^{\otimes n})) =0$ since $N^{\otimes n} \neq \one$ holds for any $n\neq 0$ and $Y$ is compact. 
Therefore it follows that $a_{j|n, m} = - \overline{a_{j|-n, m}}$ holds from $\{(U_j, a_{j|n, m}- (- \overline{a_{j|-n, m}}))\} \in H^0(\mathcal{U}, \mathcal{O}_Y(N^{\otimes n}))$.
One can easily check 
\[
\frac{d}{dw_{k, m+1}}\psi_{k, m}(w_{k, m+1}) = \left(1 + w_{k, m+1} \frac{d}{dw_{k, m+1}} \psihat_{k, m}(w_{k, m+1})\right) \cdot e^{\psihat_{k, m}(w_{k, m+1})}. 
\]
This calculation and Lemma \ref{1+e} lead that $\psi_{k, m}$ has no critical point in $V_{k, m+1}(\sigma_m - 2\eta_m)$.
Thus $\psi_{k, m}$ is a local homeomorphism on a neighborhood of the unit circle, from which it follows that $\psi_{k, m} (S^1) = S^1$ holds. \qed \\

\begin{lemma}\label{for_appendix}
The function $\psi_{k,m}$ is one-to-one on the unit circle $S^1 = \{|w_{k, m+1}| =1 \}$. 
\end{lemma}
{\proofname} Let $\Psi_{k, m}$ be a lift of $\psi_{k, m}$ which is described as below: 
\[
\Psi_{k, m}(x) = \fractwopii \log\psi_{k, m}(e^{\twopii x}) \ (x \in\mathbb{R}).  
\]
Calculating directly, one has
\[
\Psi_{k, m}(x) = x + \fractwopii \psihat_{k, m}(e^{\twopii x}). 
\]
From Lemma \ref{1+e}, 
\[
\frac{d\Psi_{k, m}(x)}{dx} = 1 + \left. \frac{d}{d\zeta}\psihat_{k, m}(e^\zeta) \right|_{\zeta = \twopii x} < 1 + \frac{1}{1+e^{\sigma_0}} < 2
\]
and
\[
\frac{d\Psi_{k, m}(x)}{dx} > 0
\]
hold. Supposing the degree of $\psi_{k, m}|_{S^1}$ is larger than 2,  it contradicts this estimate from the mean value theorem. \qed \\

{
To combine Lemma \ref{S1} and Lemma \ref{for_appendix}, we can see $\psi_{k, m} \in \Diff$. 
By considering Rouch\'{e}'s theorem, we can check the well-definedness of $\psi_{k, m}^{-1}$ to prove (3) and (4) of Proposition \ref{well-definedness}}. 
{
We will show the consequence of checking the well-definedness of the map $\psi_{k, m}^{-1}$ on $V_{k, m}(\sigma_m - 2\eta_m)$. 
In what follows, we use the notation in proof of Lemma \ref{for_appendix}. 
Let $\Psi_{k, m}$ be a lift of the map $\psi_{k, m}|_{S^1}$. 
{For a given point $w_0 \in V_{k, m}(\sigma_m -2\eta_m)$}, let $\widetilde{w}_0$ be a point in $B(\sigma_m - 2\eta_m) \coloneqq \{ | \Im\, z| < (\sigma_m - 2\eta_m)/2\pi\}$ which satisfies $e^{\twopii \widetilde{w}_0} = w_0$. 
We define a domain $D \subset \mathbb{C}$ by 
\[
D = \left\{z \in\mathbb{C} \ ; |\Im\, z| < \frac{\sigma_m -\eta_m}{2\pi}, |\Re\, z - \Re\, \widetilde{w}_0| < \frac{1}{2}\right\} 
\]
(see Figure \ref{D}). 
We define holomorphic functions $g_1, g_2$ on a neighborhood of $D$ by
\begin{gather*}
g_1(z) = z - \widetilde{w}_0 \\
g_2(z) = \Psi_{k, m}(z) - z
\end{gather*}
\begin{lemma}\label{A1}
For any $z \in \del D$, $|g_1(z)| > |g_2(z)|$ holds. 
\end{lemma}
{\proofname}
One has
\[
|g_1(z)| \geq \min \left\{ \frac{1}{2}, \ \frac{1}{2\pi}\eta_m \right\}. 
\] 
From the definition, $\{\eta_m\}$ is monotonically decreasing. 
Thus, one obtains 
\[
\frac{1}{2\pi} \eta_m < \frac{1}{2\pi} \eta_0 < \frac{1}{2}
\]
under the assumption $\eta_0 < \pi$. 
Considering the relation 
\[
\Psi_{k, m}(z) = z + \fractwopii \psihat_{k, m}(e^{\twopii z}), 
\]
one has
\[
|g_2(z)|\  {\leq } \  \frac{1}{2\pi} \sup_{z \in \del D}|\psihat_{k, m}(e^{\twopii z})| < \frac{1}{2\pi}\eta_m. 
\]
Therefore, $\displaystyle{|g_1(z)| \geq \frac{1}{2\pi}\eta_m > |g_2(z)|}$ holds. \qed \\
\begin{figure}[H]
\centering
\includegraphics[width = 10cm]{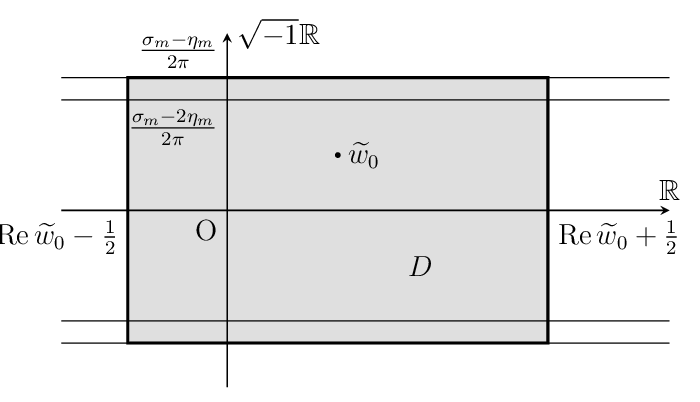}
\caption{The domain $D$ is independent of the imaginary part of $\widetilde{w}_0$. }
\label{D}
\end{figure}

{\sc Proof of (3) and (4). }
From Lemma \ref{A1}, applying Rouch\'{e}'s theorem, one finds that $\Psi_{k, m}(z) = \widetilde{w}_0$ and $z = \widetilde{w}_0$ have the same number of zeros in $D$. 
Thus the relation {{$\#\{z\in D\mid\Psi_{k, m}(z) = z\} = \#\{z\in D\mid z = \widetilde{w}_0\} =1$} } holds. 
Consequently one can define $\psi^{-1}_{k, m}$ on $V_{k, m}(\sigma_m - 2\eta_m)$.  \qed
}
%
%
\subsection{The estimate of renewed transition functions and proof of Theorem \ref{main_theorem2}}\label{new}
The renewed transition function $f_{kj, m+1}$ satisfies 
\[
\fhat_{kj, m+1}(w_{j, m+1}) = \fhat_{kj, m}(w_{j, m}) + \psihat_{j, m}(w_{j, m+1}) - \psihat_{k, m}({w_{k, m+1}}). 
\]
Combining (\ref{Schroder_eq}) and this equation, we obtain
	\begin{align}\label{g}
|\fhat_{kj, m+1}(w_{j, m+1})| &\leq |\fhat_{kj, m}(w_{j, m}) - \fhat_{kj, m}(w_{j, m+1})|  \ \ \ \ \ \ \ \\
& \ \ \ \ \ \ \ \ \ \ + |\psihat_{k, m}( e^{\alpha_{kj}}w_{j, m+1}) - \psihat_{k, m}(w_{k, m+1})|. \nonumber
	\end{align}
	\begin{claim}
If the inductive assumption $\displaystyle{\max_{j, k}||\fhat_{kj, m}||_{\sigma_m} < \delta_m}$ holds, then the renewed function $f_{kj, m+1}$ satisfies $\displaystyle{\max_{j, k}||\fhat_{kj, m+1}||_{\sigma_{m+1}} < \delta_{m+1}}$. 
	\end{claim}
{\proofname} On the first term of the right hand side of (\ref{g}), for $w_{j, m+1} \in V_{j, m+1}(\sigma_m - 4\eta_m)$, one can estimate as the following: 
	\begin{align*}
| \fhat_{kj, m}(w_{j, m}) - \fhat_{kj, m}(w_{j, m+1})| 
&= {\left| \int_{\log w_{j, m+1}}^{\log w_{j, m}} \frac{d}{d\zeta} \fhat_{kj, m}(e^\zeta) d\xi\right| }\\
&\leq \left| \log\frac{w_{j, m}}{w_{j, m+1}} \right| \cdot \sup_{|\mathrm{Re}\zeta|<\sigma_m-\eta_m} \left| \frac{d}{d\zeta}\fhat_{kj, m}(e^\zeta)\right| \\
&\leq ||\psihat_{j, m}||_{\sigma_m-4\eta_m} \cdot \sup_{|\mathrm{Re}\zeta|<\sigma_m-\eta_m} \left| \frac{d}{d\zeta}\fhat_{kj, m}(e^\zeta)\right| .
	\end{align*}
From Lemma \ref{estimate_psi}, we can easily show $\displaystyle{||\psihat_{j, m}||_{\sigma_m-4\eta_m} < \frac{C_1 \delta_m }{\eta_m^{\mu}}}$. 
By using Cauchy's integral expression, the {following} holds:
	\begin{align*}
\sup_{|\mathrm{Re}\zeta|<\sigma_m-\eta_m} \left| \frac{d}{d\zeta}\fhat_{kj, m}(e^\zeta)\right| 
&\leq \sup_{|\mathrm{Re}\zeta|<\sigma_m-\eta_m} |e^\zeta| \cdot \left| \fractwopii \int_{|\xi-e^\zeta|=\eta_m} \frac{\fhat_{kj, m}(\xi)}{(\xi-e^\zeta)^2} d\xi \right| \\
&\leq \frac{e^{\sigma_m} \delta_m}{2\pi} \int_0^{2\pi}\frac{1}{\eta_m^2}\cdot \eta_m \ d\theta \\
&=\frac{e^{\sigma_m} \delta_m}{\eta_m} < \frac{e^{\sigma_0} \delta_m}{\eta_m}.
	\end{align*}
Hence one has
\[
| \fhat_{kj, m}(w_{j, m}) - \fhat_{kj, m}(w_{j, m+1})| < \frac{C_1 \delta_m }{\eta_m^{\mu}} \cdot \frac{e^{\sigma_0} \delta_m}{\eta_m} = \frac{e^{\sigma_0}C_1 \delta_m^2}{\eta_m^{\mu+1}} .
\]
Next, for the second term of (\ref{g}), we have the estimate as below by using Lemma \ref{1+e}:
	\begin{align*}
|\psihat_{k, m}&(e^{\twopii\alpha_{kj}}w_{j, m+1}) - \psihat_{k, m}(w_{k, m+1})| 
= \left| \int_{\alpha_{kj}+\log w_{j, m+1}}^{\log w_{k, m+1}} \frac{d}{d\zeta}\left( \psihat_{k, m}(e^\zeta)\right) d\zeta \right| \\
&\leq \left| \log w_{k, m+1} - (\alpha_{kj}+\log w_{j, m+1}) \right| \cdot \sup_{|\mathrm{Re}\zeta| < \sigma_m - \eta_m} \left| \frac{d}{d\zeta}\left( \psihat_{k, m}(e^\zeta)\right) \right| \\
&< ||\fhat_{kj, m+1}||_{\sigma_m-4\eta_m} \cdot \frac{1}{1+e^{\sigma_0}}. 
	\end{align*}
Therefore, it follows that
\[
||\fhat_{kj, m+1}||_{\sigma_m - 4\eta_m} \leq \frac{e^{\sigma_0}C_1 \delta_m^2}{\eta_m^{\mu+1}} + \frac{||\fhat_{kj, m+1}||_{\sigma_m - 4\eta_m}}{1 + e^{\sigma_0}}.
\]
Solving for $||\fhat_{kj, m+1}||_{\sigma_m - 4\eta_m}$ gives
\[
\max_{j, k} ||\fhat_{kj, m+1}||_{\sigma_{m+1}} <  (1 + e^{\sigma_0})C_1 \frac{\delta_m^2}{\eta_m^{\mu+1}} =  \delta_{m+1}. 
\] \qed \\

Finally we need to prove $\delta_{m+1} < \eta_{m+1}$ and $\displaystyle{ \delta_{m+1} < \frac{\eta_{m+1}^{\mu+1}}{(1+e^{\sigma_0})C_1 \mu}}$ under the inductive assumption. 
The {latter} inequality can be proven as below:
\begin{align*}
\delta_{m+1} =   (1 + e^{\sigma_0})C_1 \frac{1}{\eta_m^{\mu +1}} \delta_m^2 &<  (1 + e^{\sigma_0})C_1 \frac{1}{\eta_m^{\mu +1}} \left( \frac{\eta_m^{\mu+1}}{(1+e^{\sigma_0})C_1\mu}\right)^2 \\
&= \frac{1}{(1+e^{\sigma_0})C_1\mu}\cdot \frac{1}{\mu}\eta_{m}^{\mu+1}\\
&= \frac{\eta_{m+1}^{\mu+1}}{(1+e^{\sigma_0})C_1 \mu}. 
\end{align*}
The {former inequality} is shown as
\begin{align*}
\delta_{m+1} &<  (1 + e^{\sigma_0})C_1 \frac{1}{\eta_m^{\mu +1}} \cdot \frac{\eta_m^{\mu+1}}{ (1 + e^{\sigma_0})C_1\mu} \cdot \delta_m \\
&= \mu^{-1}\delta_m \\
&< \mu^{-\frac{1}{\mu+1}} \eta_m \\
&= \eta_{m+1}. 
\end{align*}
{
From the assumption of Theorem \ref{main_thm}, the initial transition function $f_{kj, 0}$ satisfies

\[
\max_{j, k} ||\fhat_{kj, 0}||_{\sigma_0}= \max_{j, k} \sup_{e^{-\sigma_0}<|p|<e^{\sigma_0}} \left| \sum_{n\neq0} \, b_{kj|n, 0} \, p^n \right| < \min \left\{ \eta_0, \ \frac{\eta_0^{\mu+1}}{(1+e^{\sigma_0})C_1\mu}\right\}= \delta_0. 
\] 
Hence, we can obtain the inequality $\max_{j, k} ||\fhat_{kj, m}||_{\sigma_{m}} < \delta_m$ inductively. 
}

It is easily checked that $\delta_m$ converges to zero by $\delta_m < \eta_m$ and the limit of $\sigma_m$ is non-zero {as follows}. 

Since $\mu^{-\frac{1}{\mu+1}} < 1$ for any $\mu>1$, $\lim_{m\to\infty}\eta_m =0$ holds from the definition of $\{\eta_m\}_{m=0}^\infty$.  
Therefore we can see $\lim_{m\to\infty}\delta_m = 0$ from the inequalty $\delta_m < \eta_m$. 
The limit of $\sigma_m$ is computed as 
\[
\lim_{m\to\infty}\sigma_m = \sigma_0 - 4(\eta_0 + \eta_1 + \cdots) = \sigma_0 - \frac{4\eta_0}{1-\mu^{-\frac{1}{\mu+1}}}. 
\]
From the condition which $\eta_0$ satisfies in Theorem \ref{main_theorem2}, this limit is a non-zero constant. 
}
Therefore the main theorem follows from Theorem \ref{main_theorem2} \qed 
\section{Example}\label{example}
\begin{figure}[H]
\begin{minipage}{0.49\hsize}
\includegraphics[width = 7.5cm]{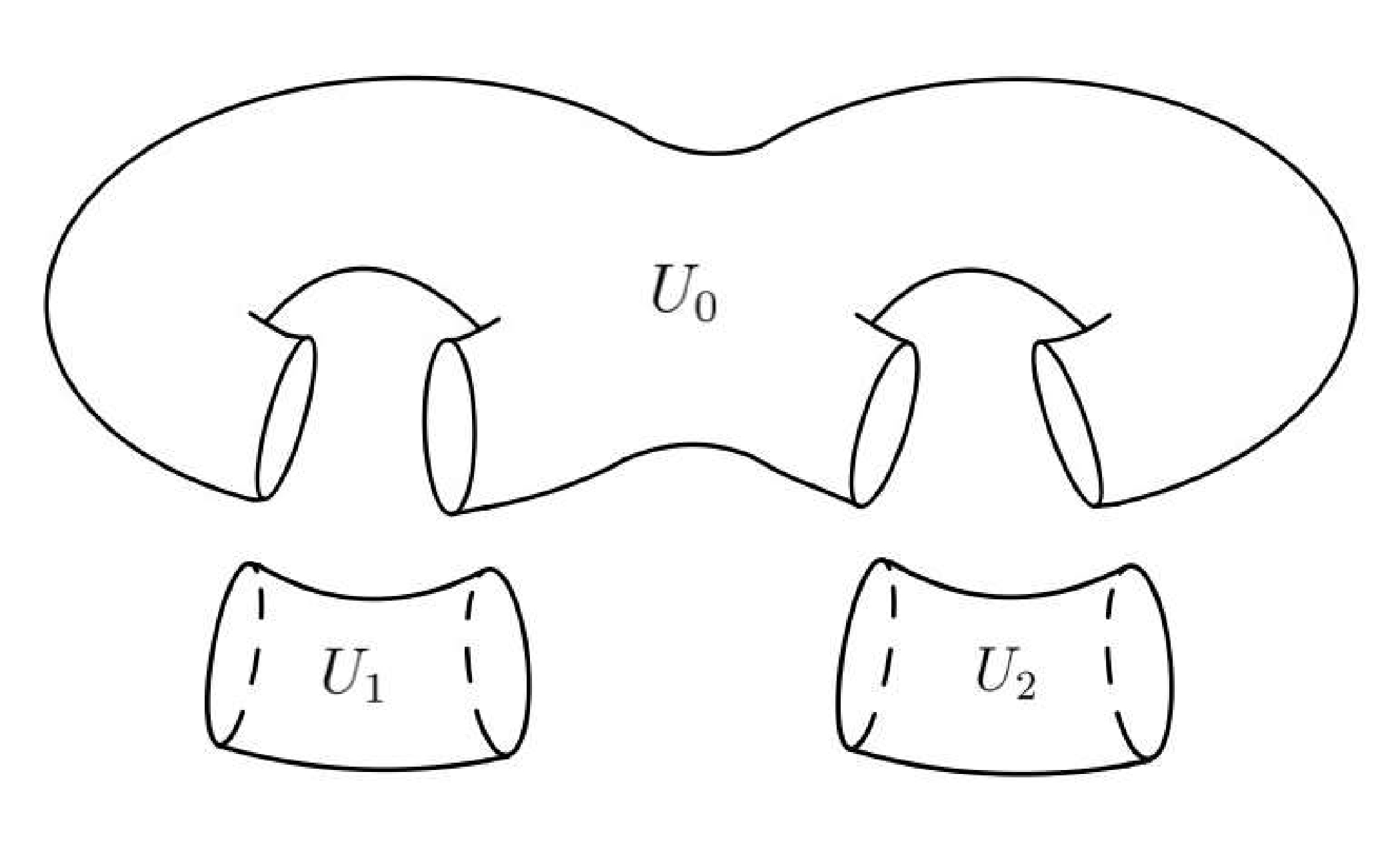}
\caption{}
\label{figure1}
\end{minipage}
\begin{minipage}{0.49\hsize}
\includegraphics[width = 7.5cm]{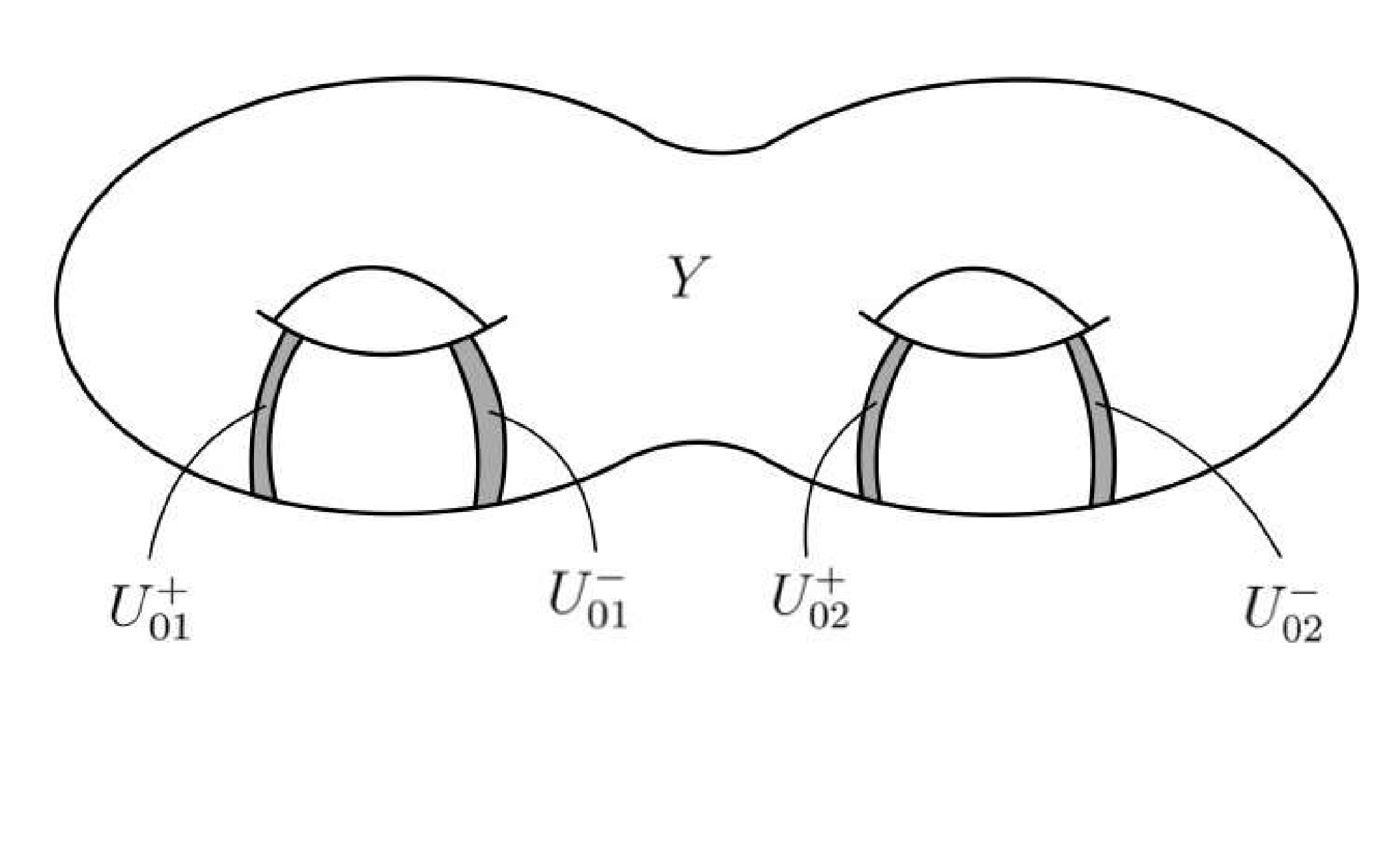}
\caption{}
\label{figure2}
\end{minipage}
\end{figure}

\begin{figure}[H]
\centering
\includegraphics[width = 7.5cm]{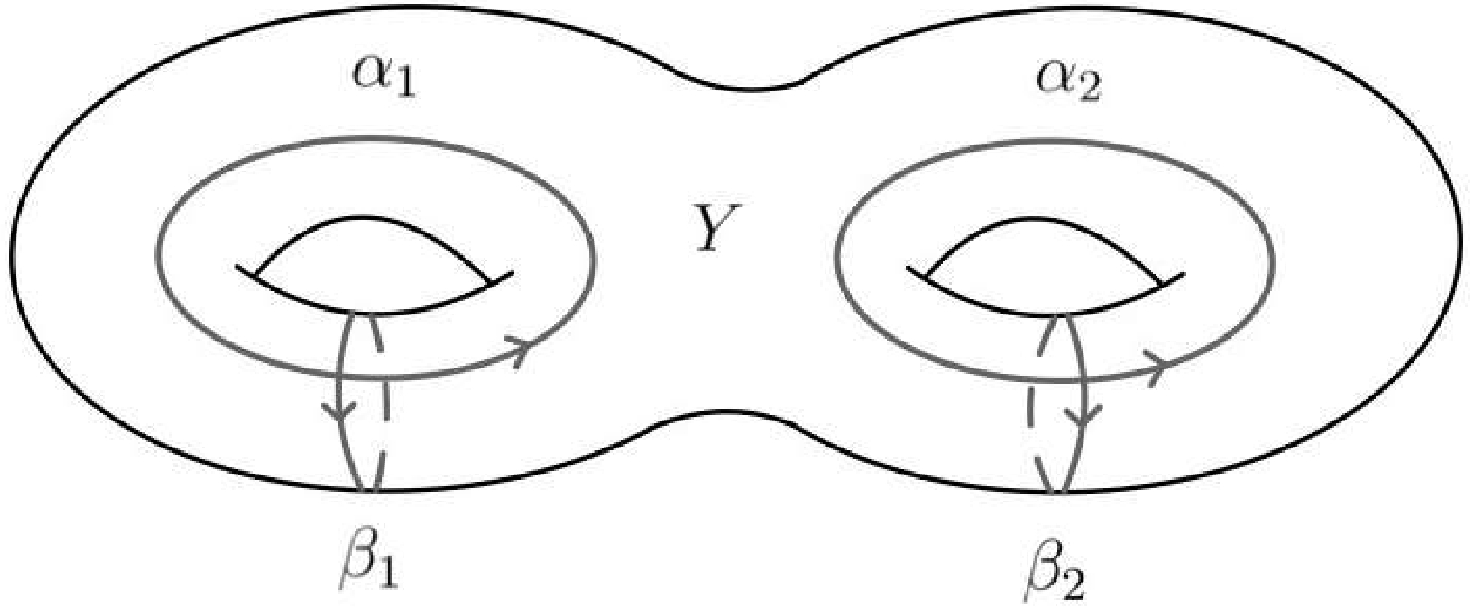}
\caption{}
\label{figure3}
\end{figure}

{
In \cite[Cor 1. ]{FK}, B.Fayad and K. Khanin showed that the familiy of commuting circle diffeomorphisms is simultaneous linearizable when rotation numbers of them satisfies the simultaneously Diophantine condition. 
}
In this section, we give a simple example and see that we get a sufficient condition for simultaneous linearization {of the pair of circle diffeomorphisms not necessarily commutative} as a consequence of the main theorem. 
Let $f_1$ and $f_2$ be elements of $\Diff$. 
For the simultaneous linearization of $f_1$ and $f_2$, we construct a Levi-flat hypersurface which has a structure of $S^1$-bundle as below. 
Let $Y$ be a compact Riemann surface with genus $2$. 
We give a finite covering $\mathcal{U} =\{U_j \}_{j = 0, 1, 2}$ of $Y$ as below (see Figure \ref{figure1}). 
The intersections are denoted by $\{U_{0j}^+, U_{0j}^-\}_{j = 1, 2}$ as Figure \ref{figure2}. \\
Define a fundamental group action $\kappa\colon\pi_1(Y, \ast) \to \Diff$ by letting $\kappa(\alpha_1) = f_1$, $\kappa(\alpha_2) = f_2$, and $\kappa(\beta_1) = \kappa(\beta_2) = \mathrm{id}_{S^1}$ for generating loops $\alpha_1, \alpha_2, \beta_1$, and $\beta_2$ (see Figure \ref{figure3}). 
By considering extending a domain of $f_1$ and $f_2$ along the unit circle, one obtains a non-singular complex surface $X$ which has a Levi-flat hypersurface $M$ constructed by suspension of $\kappa$. 
Let $\pi \colon M \to Y$ be the projection and $P$ be a holomorphic submersion as in \S1. 
For the good system of local functions $\{(V_j, (z_j, w_j))\}$ of $\sigma > 0$ over $\{Y, \mathcal{U}, \kappa, X, P\}$, we define the transversal transition on each $V_{01}^+, V_{01}^-, V_{02}^+$, and $V_{02}^-$, where $\{V_{0j}^\epsilon\}_{\epsilon=+, -, \ j= 1,2}$ are the connected components of the intersection {$V_{0j}= V_0 \cap V_j = P^{-1}(U_0) \cap P^{-1}(U_j)\ (j = 1, 2)$}.  
Denote by $w_j = f_{j0}^+(w_0) = f_j(w_0)$ the transition on $V_{0j}^+$ and by $w_j = f_{j0}^-(w_0) = w_0$ the transition $V_{0j}^-$ for each $j = 1, 2$. 
Then, the Laurent expansions for transition functions are written as below: 
\begin{align*}
\log\frac{f_{j0}^+(w_0)}{w_0} &= \log\frac{f_j(w_0)}{w_0} \equiv \alpha_{j0}^+ + \sum_{n\neq0}b_{j0|n}^+w_0^n, \\
\log\frac{f_{j0}^-(w_0)}{w_0} &= \log\frac{\mathrm{id}(w_0)}{w_0}\equiv \alpha_{j0}^- \equiv 0. 
\end{align*}
In this situation, since $U_0 \cap U_1 \cap U_2 = \emptyset$, any \v{C}ech 1-cochain over $Y$ satisfies 1-cocycle condition. 
Therefore one obtains 
\[
N = [\alpha(\{(V_j, (z_j, w_j))\})] = [\{(U_{0j}^+, e^{\alpha_{j0}^+}), (U_{0j}^-, 1)\}_{j = 1, 2}]\in \check{H}^1(\mathcal{U}, \U) 
\] and 
\[
[b_n(\{(V_j, (z_j, w_j))\})] = [\{(U_{0j}^+, b_{j0|n}^+), (U_{0j}^-, 0)\}_{j = 1, 2}]. 
\]
\begin{proposition}
Let $f_1$, $f_2$, $Y$, $\mathcal{U}$, $\kappa$, $X$, and $P$ be as above and $\{(V_j, (z_j, w_j)) \}_{j = 0, 1, 2}$ be the good system of local functions of width $\sigma$. 
{ If the following conditions $(i)$, $(ii)$ and $(iii)$ hold}, then $f_1$ and $f_2$ are simultaneous linearizable. \\
$(i)$ The unitary flat line budle $N$ over $Y$ satisfies $(C_0, \mu, K)$-Diophantine condition with the constant $K = K(Y, \mathcal{U})$ obtained by Theorem \ref{U}, where $C_0 > 0$ and $\mu >1$ is a constant determined only by $Y$ and $\mathcal{U}$. \\
$(ii)$ {
For non-zero order coefficients $b_{j0|n}^+$ associated to the transversal transition function $f_j$ of $\{(V_{0j}^\epsilon(z_j, w_j))\}_{\epsilon=+, -, j=1, 2}$, there exists a constant $\eta_0 \in (0, \min\{\pi, (1-\mu^{-\frac{1}{\mu+1}})\frac{\sigma}{4}\})$ such that
\[
\max_{j =1, 2}\sup_{e^{-\sigma} < |p| < e^\sigma}\left| \sum_{n \neq 0} b_{j0|n}^+ p^n \right| <\min\left\{ \eta_0, \ \frac{\eta_0^{\mu+1}}{(1+e^{\sigma})C_1\mu}\right\}
\]
holds, where the constant $C_1$ is obtained by Lemma \ref{estimate_psi}. \\
}
 $(iii)$ The 1-cohomology group $[b_n(\{(V'_j, (z_j, w'_j))\})]$ is cohomologous to $0$ for any good system of local functions $\{(V'_j, (z_j, w'_j))\}$ which satisfies $[\alpha(\{(V'_j, (z_j, w'_j))\})] = N$. \\
\end{proposition}
{\proofname}
From Theorem \ref{main_thm}, the system $\{Y, \mathcal{U}, \kappa, X, P, \{(V_j, (z_j, w_j)) \}, \sigma \}$ is linearizable. 
Then, there exist the functions of retaking coordinates $\psi_j (u_j) = w_j \ (j = 0, 1, 2)$ and the good system of local function $\{(V'_j, (z_j, u_j))\}_{j = 0, 1, 2}$ whose transition is the linear map $u_j = e^{\alpha_{j0}^\epsilon} u_0$ for each $\epsilon = +, -$ and $j=1, 2$. 
One obtains the following relations:
\[
(\psi_j^{-1} \circ f_{j0}^\epsilon \circ \psi_0)(u_0) = e^{\alpha_{j0}^\epsilon} u_0. 
\]
Therefore we obtain the {following} for each $j=1, 2$:
\begin{gather*}
(\psi_j^{-1} \circ f_{j} \circ \psi_0)(u_0) = e^{\alpha_{j0}^\epsilon} u_0 \\ 
(\psi_j^{-1} \circ \psi_0)(u_0) = u_0.
\end{gather*}
Hence, one can check $(\psi_0^{-1} \circ f_{j} \circ \psi_0) (u_0) = e^{\alpha_{j0}^+} u_0$ for each $j = 1, 2$. 
That is simultaneous linearization of $f_{1}$ and $f_{2}$. 
Since $f_{j}$ is conjugated to the rotation $w \mapsto e^{\alpha_{j0}^+}w$, $\alpha_{j0}^+ \equiv \twopii\rho(f_{j}) \modtwopiiZ$ holds, where $\rho(f_{j})$ is the rotation number of  $f_{j}$. \qed \\

\section{Discussion}\label{discussion}
In this paper, we considered a function $f \in \Diff$ which has the Laurent expansion
\[
\log \frac{f(w)}{w} = \alpha + \sum_{n \neq 0}b_n w^n. 
\]
Comparing the condition of main theorem with the assumption in Arnol'd's linearization theorem, we should check a relation of the constant term of Laurent expansion as above and the rotation number of $f$. 
As is seen the example in \S\ref{example}, sometimes $\alpha$ turns out to be equal to $\twopii\rho(f)$ modulo $2\pi\sqrt{-1}\mathbb{Z}$. 
\begin{question}
Assume that $\alpha$ is a constant term Laurent expansion of $\displaystyle{\log\frac{f(w)}{w}}$. 
Then, will $e^{\alpha} = e^{\twopii \rho(f)}$ always hold?  
\end{question}
\vskip3mm
{\bf Acknowledgment. } 
The author would like to thank Takayuki Koike for fruitful comments. 
This work was partly supported by Osaka Central Advanced Mathematical Institute: MEXT Joint Usage/Research Center on Mathematics and Theoretical Physics JPMXP0619217849.

\end{document}